\documentclass[11pt]{article}
\usepackage{setspace}
\usepackage{amsfonts}
\usepackage{amssymb}         
\usepackage{amsmath}    
\usepackage{multirow}
\usepackage{graphicx}
\usepackage{algorithm}
\usepackage{algorithmic}

\newcommand{\Eb}{{\mathbf{E}}}

\newtheorem{lemma}{Lemma}

\newtheorem{proposition}{Proposition}
\newtheorem{remark}{Remark}

\numberwithin{equation}{section}
\newtheorem{thm}{Theorem}[section]
\begin{document}
\setcounter{footnote}{0}
\title{An Algorithm to Estimate a Nonuniform Convergence Bound in the Central Limit Theorem} 
\author{\textsc{Vladimir Nikulin\thanks{Email: vnikulin.uq@gmail.com}}    \\ 
Department of Mathematics, University of Queensland \\
Brisbane, Australia
}
\date{}
\maketitle
\thispagestyle{empty}
\begin{abstract}
A nonuniform version of the Berry-Esseen bound was proved.
The most important feature of the new bound is a
 monotonically  decreasing function $C(|t|)$ instead of the universal constant
$C=29.1174$:
$C(|t|)<C$ if $|t|\geq 3.2$ and $C(|t|) \rightarrow 1+e$ if $|t| \rightarrow \infty$
 where $t$ is a coordinate of the point.

 The function $C(|t|)$ has very complex analytical expression based on indicator functions.
A general algorithm was developed in order to estimate values of $C(|t|)$ for
an arbitrary $t$.
\end{abstract}
\section{Introduction}  \label{intro}
 Much of the importance of the central limit theorem follows  from its proven
 adaptability and utility in many areas of mathematics, probability theory
 and statistics \cite{Hal82}.
 The origin of the term "central" is not clear. Some authors allude to the central
 role of the theorem in probability theory, while others refer to the fact that
 the theorem concerns a measure of central tendency, namely the mean of a
 normalized sum of random variables.

Let $X,X_1,X_2, \cdots$ be a sequence of independent and identically distributed random
variables with
\begin{equation} \label{eq:beress}
 \Eb X =0, \Eb X^2 =1, \Eb |X|^3 = \rho < \infty.
\end{equation}

The classical Berry-Esseen Theorem (\cite{Ber41} and \cite{Ess45}) states that
\begin{equation} \label{eq:berry}
   H_n(t) = \frac{\sqrt{n}}{\rho} | F_n(t) - \Phi(t) | \leq A
   \hspace{0.1in} \forall n \geq 1
\end{equation}
where $A$ is an absolute constant, $F_n$ is the distribution function of
$\frac{1}{\sqrt{n}} \sum_{j=1}^n X_j$, and
$\Phi$ is a standard normal distribution function.

\cite{Nag65} and \cite{Nag79}  established a nonuniform structure of the Berry-Esseen bound:
\begin{equation} \label{eq:nagev}
   H_n(t) \leq \frac{C}{1+|t|^3}  \hspace{0.1in} \forall n \geq 1.
\end{equation}

\cite{Nik92} and \cite{NikPad98} demonstrated that in fact the constant
$C$ in the upper bound (\ref{eq:nagev}) may be replaced by a decreasing 
function of $|t| \geq 1$:
\begin{equation} \label{eq:nikul}
  H_n(t) \leq C(|t|) \cdot |t|^{-3} \hspace{0.1in} \forall n \geq 1.
\end{equation}

The best known values of the above constants are as follow
\begin{equation} \label{eq:shig}
   A \leq 0.7655;
\end{equation}
\begin{equation} \label{eq:mich}
   C \leq 29.1174 \approx 0.7655 \cdot \left( 1+ \left(\frac{10}{3} \right)^3 \right)
\end{equation}
where constants (\ref{eq:shig}) and (\ref{eq:mich}) were proved by \cite{Shi82} and
 \cite{Mic81}.

\section{Main Results}  \label{sec:main}
In this section we formulate the algorithm in order to estimate
values of the function $C(|t|)$ for an arbitrary value of the argument.
Without loss of generality we consider the case $t>0$ only,
because the estimation procedure is symmetrical against the point of origin.

\begin{thm} \label{th:mainth} Under conditions of the Berry-Esseen Theorem (\ref{eq:beress})
Tables \ref{tb:table1} and \ref{tb:table2} represent values for the upper bound
\[ \underset{x \geq t}{\sup} \{ \frac{\sqrt{n}x^3}{\rho} |F_n(x) - \Phi(x)| \} \leq C(t)
  = \max{\{B_T(t), B_C(t)\}}  \]
where functions in the right part of the above equation are defined in (\ref{eq:frst2}) and
(\ref{eq:form7} - \ref{eq:form7x}).
\end{thm}

Very briefly, the proof
 structure of the Theorem~\ref{th:mainth} may be formulated as follows.

The target is to estimate
$\underset{t \geq t_0}{\sup} \hspace{0.02in} \{ t^3 H_n(t)\},$ and we can exploit
exponential rate of decline to zero of the standard normal distribution function
$1-\Phi(t)$ if $t \rightarrow \infty.$
Respectively, we will split the task into 2 parts:
\renewcommand{\theenumi}{A\arabic{enumi}}
\renewcommand{\labelenumi}{\theenumi)}
\begin{enumerate}
  \item \texttt{tail}: $[\psi(n, t_0), \infty[$, see Lemma~\ref{th:first};

  \item \texttt{center}: $[t_0, \psi(n, t_0)]$, see Lemma~\ref{th:second}
\end{enumerate}
where $\psi(n, t_0)$ is an increasing function of both arguments $t_0$ and $n$
(to be defined in (\ref{eq:frst1})).

Firstly,
we construct the upper bound $B_T(t)$ for
$\frac{\sqrt{n}}{\rho} t^3 (1-F_n(t)), t \geq \psi(n, t_0).$
This bound has an essential property: $B_T(t) \geq \frac{\sqrt{n}}{\rho} t^3 (1-\Phi(t))
\hspace{0.02in} \forall t \in [\psi(n, t_0), \infty[.$

The second step is a much more complex. Using truncation method and results of the
Lemma~\ref{th:trunk} we construct the upper bound
$B_C(t)$ for $t^3 H_n(t), t_0 \leq t \leq \psi(n, t_0).$

Subject to the special conditions, both bounds $B_T(t)$ and $B_C(t)$ represent 
decreasing functions of $t$ and independent on $n$.

Finally, the required value $C(t_0)$ will be computed as a maximum
of the ``tail'' and ``center'' bounds.

The Table~\ref{tb:table1} demonstrates 1) improvement of the nonuniform bound (\ref{eq:nagev})
with the constant (\ref{eq:mich}) if $|t| \geq 3.2$; 2) improvement of the uniform bound
(\ref{eq:berry}) with the constant (\ref{eq:shig}) if $|t| \geq 3.3$.

\begin{table}
\caption{Values of the upper bound $C(t)$ where $\tau$ and $b$ are
components of the truncation parameter $h$ to be defined in the Section~\ref{sec:main}.}
   \label{tb:table1}
\begin{center}
 \begin{tabular}{c|cc|c|cc}
   t & $\tau$ & b & C(t) & Bound: (\ref{eq:nikul}) & Bound: (\ref{eq:nagev}) \\
\hline
\noalign{\smallskip}
 3.18 & 0.4553 & 1.9690 & 28.4057 & 0.88333358 & 0.87823442 \\
 3.19 & 0.4570 & 1.9670 & 28.3187 & 0.87237015 & 0.87024709 \\
 3.20 & 0.4587 & 1.9650 & 28.2363 & \textbf{0.8617025} & 0.86235485 \\
 3.21 & 0.4604 & 1.9637 & 28.1563 & 0.85125797 & 0.85455633 \\
 3.22 & 0.4601 & 1.9617 & 28.0809 & 0.84109132 & 0.84685015 \\
 3.23 & 0.4588 & 1.9597 & 28.0052 & 0.83105872 & 0.83923498 \\
 3.24 & 0.4584 & 1.9577 & 27.9293 & 0.82115399 & 0.83170951 \\
 3.25 & 0.4581 & 1.9557 & 27.8532 & 0.81138124 & 0.82427245 \\
 3.26 & 0.4577 & 1.9547 & 27.7743 & 0.80166121 & 0.81692251 \\
 3.27 & 0.4563 & 1.9527 & 27.6980 & 0.79214601 & 0.80965846 \\
 3.28 & 0.4559 & 1.9507 & 27.6215 & 0.78275383 & 0.80247906 \\
 3.29 & 0.4555 & 1.9487 & 27.5448 & 0.77348593 & 0.79538311 \\
 3.30 & 0.4551 & 1.9467 & 27.4681 & \textbf{0.7643403} & 0.78836942 \\
 3.40 & 0.4506 & 1.9284 & 26.6933 & 0.67915056 & 0.72250902 \\
 3.50 & 0.4461 & 1.9097 & 25.9186 & 0.60451529 & 0.66370384 \\
 3.60 & 0.4439 & 1.8907 & 25.1491 & 0.53903300 & 0.61104599 \\
 3.70 & 0.4419 & 1.8717 & 24.3886 & 0.48148354 & 0.56376222 \\
 3.80 & 0.4402 & 1.8517 & 23.6406 & 0.43083173 & 0.52119151 \\
 3.90 & 0.4399 & 1.8327 & 22.9052 & 0.38613551 & 0.48276685 \\
 4.00 & 0.4400 & 1.8137 & 22.1853 & 0.34664525 & 0.44800024 \\
 4.10 & 0.4405 & 1.7939 & 21.4825 & 0.31169807 & 0.41647019 \\
 4.20 & 0.4406 & 1.7749 & 20.7969 & 0.28070572 & 0.38781182 \\
 4.30 & 0.4412 & 1.7559 & 20.1302 & 0.25318851 & 0.36170787 \\
 4.40 & 0.4434 & 1.7359 & 19.4829 & 0.22871601 & 0.33788194 \\
 4.50 & 0.4441 & 1.7179 & 18.8552 & 0.20691545 & 0.31609246 \\
 4.60 & 0.4465 & 1.6989 & 18.2478 & 0.18747242 & 0.29612776 \\
 4.70 & 0.4485 & 1.6809 & 17.6609 & 0.17010635 & 0.27780183 \\
 4.80 & 0.4502 & 1.6619 & 17.0947 & 0.15457488 & 0.26095081 \\
 4.90 & 0.4525 & 1.6449 & 16.5490 & 0.14066430 & 0.24543000 \\
 5.00 & 0.4556 & 1.6269 & 16.0240 & 0.12819173 & 0.23111131 \\
 \end{tabular}
\end{center}
\end{table}

All values in the Tables \ref{tb:table1} and \ref{tb:table2} were computed
using the following Algorithm where $\tau(t)$ and $b(t)$ are an important components
of the truncation parameter $h(t)$ to be defined in (\ref{eq:frst1}) and (\ref{eq:frst1c}).
\newline
\vspace{0.2in}
\newline
\textbf{Algorithm 1.} \small{(for computation of the values of the upper bound $C(t)$
 in the Tables~\ref{tb:table1} and ~\ref{tb:table2})}
\vspace{0.05in}
\renewcommand{\theenumi}{\arabic{enumi}}
\renewcommand{\labelenumi}{\theenumi:}
\begin{enumerate}
\item Enter value of the argument $t \geq 3.18$;
\item compute
  \[ \overline{\tau}(t) := \min{\{0.5(1+\sqrt{1-\frac{10}{t^2}}), 1-\frac{\sqrt{3}}{t} \}} \]
 (this step follows from conditions (\ref{eq:gcond1}) and (\ref{eq:gcond7}));

\item compute \[ \overline{b}(t) := \sqrt[3]{\frac{30}{1+e}} \]
  (this step follows from the structure of the ``tail'' bound (\ref{eq:frst2}));

\item compute
  \[ \underline{\tau}(t) := \max{\{ \tau_1, 0.5(1-\sqrt{1-\frac{10}{t^2}}) \}} \]
  where
  \[ \tau_1 = \frac{1+\sqrt{1-4p}}{2} \hspace{0.1in} \mbox{if}
  \hspace{0.1in}
  4p \leq 1, p=\frac{2(b(t)-1)}{b^2(t)} - \frac{1}{t^2}, \]
  \[ \mbox{alternatively,} \hspace{0.1in} \tau_1 = -\infty \]
  (this step follows from conditions (\ref{eq:gcond6}) and (\ref{eq:gcond1}));

\item compute \[ \underline{b}(t) := \frac{2t}{t+\sqrt{t^2-6}} \]

 (this step follows from condition (\ref{eq:gcond5}));
\end{enumerate}

\renewcommand{\theenumi}{\arabic{enumi}}
\renewcommand{\labelenumi}{\theenumi:}
\begin{enumerate}
\item[6:] enter parameters $\tau(t)$ and $b(t)$ within the ranges:
  $\underline{\tau}(t) \leq \tau(t) \leq \overline{\tau}(t)$ and
  $\underline{b}(t) \leq b(t) \leq \overline{b}(t);$

\item[7:] check conditions (\ref{eq:term7}), (\ref{eq:cond4}) and (\ref{eq:cond6});

\item[8:] compute upper bounds for the parameters $\gamma(t), \mu(t)$ and
$\beta(t)$ according to
 (\ref{eq:cond2}), (\ref{eq:form2}) and  (\ref{eq:form6b});

\item[9:] check conditions (\ref{eq:form4}) and (\ref{eq:form4ba} - \ref{eq:form4c});

\item[10:] compute low and upper bounds $ \underline{m}_2(t)$ and $\overline{m}_2(t) $
  for $m_2(t)$ according to (\ref{eq:form6c}) using the ``center'' condition (\ref{eq:cond2});

\item[11:] compute low bound for $\delta(t)$ according to (\ref{eq:form1}) and
 $\alpha_k(t), k=0..3, \Delta(t)$ according to (\ref{eq:alfa}) and (\ref{eq:alfa1});

\item[12:] compute $\eta(t) := |\mu(t)| \sqrt{\overline{m}_2(t)}(3\sqrt{\overline{m}_2(t)}
  +|\mu(t)| \sqrt{\beta(t)})$ using upper bounds (\ref{eq:cond2}), (\ref{eq:form2})
  and (\ref{eq:form6c}).

\item[13:] Finally, $C(t) = \max{\{B_T(t), B_C(t) \}}$
where the bound $B_T(t)$ is defined in  (\ref{eq:frst2}), and the bound $B_C(t)$ is defined in
(\ref{eq:form7} - \ref{eq:form7x}).
\end{enumerate}
\begin{remark} The pair of parameters $\left(\underline{\tau}(t), \overline{b}(t)\right)$ will
pass all required conditions if $t \geq 3.18.$ In order to optimize selection of the
parameters we can consider all possible values from the intervals
$\left[\underline{\tau}(t), \overline{\tau}(t) \right]$ and
$\left[\underline{b}(t), \overline{b}(t) \right]$ with ordered steps.
A Pentium 4, 2.8GHz, 512MB RAM, computer was used for the computations which were conducted according to the
special program written in C. The total computation time for all values in the
Tables~\ref{tb:table1} and ~\ref{tb:table2} with both steps equal to $0.001$ was less than 1 min.
\end{remark}

\begin{table}
\caption{Further values of the upper bound $C(t)$.}
           \label{tb:table2}
\vspace{0.1in}
\begin{center}
 \begin{tabular}{c|cc|c|cc}
   t & $\tau$ & b & C(t) & Bound: (\ref{eq:nikul}) & Bound: (\ref{eq:nagev}) \\
\cline{1-6}
\hline
\noalign{\smallskip}
 6.00 & 0.4843 & 1.4696 & 11.8046 & 0.05465073 & 0.13419355 \\
 7.00 & 0.5166 & 1.3450 & 9.0590 & 0.02641108 & 0.08465116 \\
 8.00 & 0.5475 & 1.2486 & 7.2512 & 0.01416244 & 0.05676413 \\
 9.00 & 0.5765 & 1.1749 & 6.0329 & 0.00827556 & 0.03989041 \\
10.00 & 0.6298 & 1.1555 & 5.7370 & 0.00573698 & 0.02909091 \\
11.00 & 0.6625 & 1.1461 & 5.5971 & 0.00420522 & 0.02186186 \\
12.00 & 0.6867 & 1.1381 & 5.4808 & 0.00317173 & 0.01684211 \\
13.00 & 0.7078 & 1.1311 & 5.3802 & 0.00244890 & 0.01324841 \\
14.00 & 0.7253 & 1.1251 & 5.2951 & 0.00192969 & 0.01060838 \\
15.00 & 0.7405 & 1.1191 & 5.2108 & 0.00154394 & 0.00862559 \\
16.00 & 0.7537 & 1.1141 & 5.1413 & 0.00125519 & 0.00710764 \\
17.00 & 0.7661 & 1.1091 & 5.0724 & 0.00103244 & 0.00592593 \\
18.00 & 0.7768 & 1.1051 & 5.0177 & 0.00086037 & 0.00499229 \\
19.00 & 0.7868 & 1.1011 & 4.9634 & 0.00072363 & 0.00424490 \\
20.00 & 0.7954 & 1.0971 & 4.9095 & 0.00061369 & 0.00363955 \\
30.00 & 0.8543 & 1.0709 & 4.5661 & 0.00016911 & 0.00107848 \\
40.00 & 0.8857 & 1.0568 & 4.3888 & 0.00006858 & 0.00045499 \\
50.00 & 0.9054 & 1.0475 & 4.2732 & 0.00003419 & 0.00023296 \\
60.00 & 0.9191 & 1.0400 & 4.1827 & 0.00001936 & 0.00013481 \\
70.00 & 0.9293 & 1.0351 & 4.1237 & 0.00001202 & 0.00008490 \\
80.00 & 0.9373 & 1.0318 & 4.0843 & 0.00000798 & 0.00005687 \\
90.00 & 0.9428 & 1.0291 & 4.0527 & 0.00000556 & 0.00003995 \\
100.00 & 0.9477 & 1.0263 & 4.0200 & 0.00000402 & 0.00002912 \\
\noalign{\smallskip}
\hline
\noalign{\smallskip}
$\infty$ & 1- & 1+ & \textbf{3.7183} & 0+ & 0+ \\
 \end{tabular}
\end{center}
\end{table}

\section{Proofs}  \label{sec:main}
The proposed method is based on the following truncation
\begin{equation*}
    Y := \begin{cases}
    X & \text{if} \hspace{0.1in} | X | \leq h \\
    0  &  \text{otherwise} \\
    \end{cases}
\end{equation*}
where $h>0$ is a truncation parameter, and may be regarded as an extension of
\cite{Mic76} and \cite{Mic81}.

We will denote by $F$ and $Q$ distribution
functions of random variables $X$ and $Y$.

\setcounter{thm}{0}
\begin{lemma} \label{th:markov} (\texttt{Markov Inequality})
Suppose that $\ell$ is an arbitrary non-decreasing and
non-negative function. Then, for any $h > 0:$
\begin{equation} \label{eq:markdef}
  \mathbb{P}(|X| \geq h) \leq \frac{\Eb \ell(|X|)}{\ell(h)}, \ell(h) > 0.
\end{equation}
\end{lemma}

\begin{lemma} \label{th:trunk} (\texttt{Truncation})
The following upper bounds are valid for an arbitrary parameters
$s \geq 0$ and $h > 0$:
\begin{subequations}
\begin{align}
  \label{eq:pow0}
  \beta := \Eb \exp{\{sY\}} \leq 1+\frac{s^2}{2}+\frac{\rho}{h^3}
  \left( \exp{\{sh\}} - 1 \right);  \\
  \label{eq:pow1}
  m_1 := \Eb Y \exp{\{sY \}} \leq s+\frac{\rho}{h^2} \exp{\{sh\}};  \\
\label{eq:pow2}
  m_2 := \Eb Y^2 \exp{\{sY \}} \leq 1+\frac{\rho}{h} \exp{\{sh\}}; \\
\label{eq:pow3}
  m_3 := \Eb | Y |^3 \exp{\{sY \}} \leq \rho \exp{\{sh\}};  \\
\label{eq:pow4}
  \Eb | Y | \exp{\{sY \}} \leq \sqrt{\beta m_2}.
\end{align}
\end{subequations}
\end{lemma}

\textit{Proof:} Proofs of (\ref{eq:pow0}), (\ref{eq:pow1}), (\ref{eq:pow2}) and
(\ref{eq:pow3}) are similar and based on the Taylor representation
\begin{equation} \label{eq:taylor}
    \exp{\{sY\}} =\sum_{i=0}^{\infty} \frac{(sY)^i}{i!}.
\end{equation}

The following relations are valid according to (\ref{eq:markdef})
\begin{subequations}
\begin{align}
  \label{eq:step1}
  \Eb Y \leq \Eb X + \int_{|X| \geq h} |X| F(dX) \leq \frac{\rho}{h^2};    \\
  \label{eq:step2}
  \Eb Y^2 \leq \Eb X^2 = 1 \leq 1 + \frac{\rho}{h};  \\
\label{eq:step3}
  \Eb Y^i \leq \Eb |X|^3 h^{i-3} \leq \rho h^{i-3}, i \geq 3.
\end{align}
\end{subequations}

Combining (\ref{eq:taylor}) with (\ref{eq:step1}), (\ref{eq:step2}) and
(\ref{eq:step3}) we will obtain the bounds (\ref{eq:pow0}-\ref{eq:pow3}).

The proof of (\ref{eq:pow4}) is based on the definitions (\ref{eq:pow0}) and
(\ref{eq:pow2}) and follows from H\"older's inequality. $\blacksquare$

In order to simplify notations we will omit dependence between parameters
$\tau, b, h, c, s$ and the coordinate of the point $t$.

\begin{lemma} (\texttt{Tail approximation}) \label{th:first}
Suppose that
\begin{equation} \label{eq:frst}
  \frac{t\sqrt{n}}{h}-\frac{sn}{2h}-c \geq 0 \hspace{0.1in}
  \text{or} \hspace{0.1in} t^2 \geq \varphi_{n,t}(a,b,c)
\end{equation}
where $a>0, t >1, b>c \geq 1$ and
\begin{equation} \label{eq:frst1}
  \psi(n,t) := \varphi_{n,t}(a,b,c)=\frac{b^2}{2(b-c)} \log{\frac{\sqrt{n}t^3}{\rho a}}, \hspace{0.04in}
  h=\frac{\sqrt{n}t}{b}, \hspace{0.04in}
  s=\frac{1}{h} \log{\frac{\sqrt{n}t^3}{\rho a}}>0.
\end{equation}
Then,
\begin{equation} \label{eq:frst2}
  t^3 \cdot H_n(t) \leq  B_T(t) = b^3(1+e).
\end{equation}
\end{lemma}
\textit{Proof:} According to (\ref{eq:markdef}),
$$ 1-F_n(t)  \leq 1-Q_n(t)+ 1- (1-\mathbb{P}(|X|>h))^n $$
\begin{equation} \label{eq:markov}
   \leq 1-Q_n(t)+n\mathbb{P}(|X|>h) \leq \beta^n
   \exp{\{-st\sqrt{n}\}}+\frac{b^3\rho}{\sqrt{n}t^3}
\end{equation}
where $Q_n$ is a distribution function of $\frac{1}{\sqrt{n}} \sum_{j=1}^n Y_j$.

According to (\ref{eq:pow0}),
\begin{equation} \label{eq:mark}
   \beta = \Eb \exp{\{sY\}} \leq
   1 + \frac{s^2}{2} + \frac{b^3}{an} \leq \exp{\{ \frac{s^2}{2}+\frac{b^3}{an} \}}.
\end{equation}

The following relations are valid as a consequence of the condition (\ref{eq:frst})
$$ 1-Q_n(t) \leq \beta^n \exp{\{-st\sqrt{n}\}} $$
\begin{equation} \label{eq:cond1}
\leq \left[ \frac{a\rho}{\sqrt{n}t^3} \right]^c
   \exp{\{-sh \left( \frac{t\sqrt{n}}{h}-\frac{sn}{2h}-c \right) + \frac{b^3}{a} \}} \leq
  \left[ \frac{a \rho}{\sqrt{n}t^3} \right]^c \exp{\{\frac{b^3}{a}\}}.
\end{equation}

Combining (\ref{eq:markov}) and (\ref{eq:cond1}) we obtain
\begin{equation} \label{eq:fappr}
  1-F_n(t) \leq \frac{\rho b^3}{\sqrt{n} t^3} +
  \left[ \frac{a \rho}{\sqrt{n}t^3} \right]^c \exp{\{\frac{b^3}{a}\}}.
\end{equation}

On the other hand,
$$1 - \Phi(t) \leq \frac{1}{t\sqrt{2\pi}}
 \int_{t}^{\infty} v \exp{\{-\frac{v^2}{2}\}} dv =
 \frac{1}{t\sqrt{2\pi}} \exp{\{-\frac{t^2}{2}\}} \leq
 \frac{1}{\sqrt{2\pi}} \left[\frac{a\rho}{\sqrt{n}t^3} \right]^c $$
where the last formula was obtained using condition (\ref{eq:frst}), $t \geq 1$, and
$\frac{b^2}{4(b-c)} \geq c$ if $b \geq c$.

As far as above estimator of $1 - \Phi(t)$ is smaller comparing with (\ref{eq:fappr})
we can ignore it:
\begin{equation} \label{eq:happr}
  H_n(t) \leq \frac{b^3}{t^3} + \frac{a}{t^3} 
  \left[ \frac{a \rho}{\sqrt{n}t^3} \right]^{c-1} \exp{\{\frac{b^3}{a}\}}.
\end{equation}

Maximizing (\ref{eq:fappr}) as a function of $a$ we find $a=\frac{b^3}{c}$ and
$$ \inf_{a>0} \{ a^c \exp{ \{ \frac{b^3}{a} \} \}} = \left[ \frac{b^3e}{c} \right]^c.$$

Taking into account that $\frac{a\rho}{\sqrt{n}t^3} <1$ we conclude that the
upper bound of $H_n$ is a decreasing function of $c$. We
obtain required result if $c=1$.  $\blacksquare$

\subsection{Center approximation} Suppose that
$t_0^2 \leq t^2 \leq \varphi_{n,t}(a,b,c)$ where $a = \frac{b^3}{c}$ or
\begin{equation} \label{eq:cond2}
  \frac{\rho}{\sqrt{n}} \leq c \left( \frac{t}{b} \right)^3
  \exp{\{2(c-b)\left(\frac{t}{b}\right)^2 \}} = \gamma(t).
\end{equation}

It will be more convenient for us to redefine here some of the variables of the
Lemma~\ref{th:first}:
\begin{equation} \label{eq:frst1c}
  h=\tau\sqrt{n} t,\hspace{0.02in} r=t(1-\tau), \hspace{0.02in} s=\frac{r}{\sqrt{n}},
  \hspace{0.02in} \varepsilon=\frac{\rho}{\tau^3\sqrt{n}t^3} \exp{ \{ \tau (1-\tau) t^2 \}},
  \hspace{0.02in} 0 < \tau < 1;
\end{equation}
$$G_n(x)=\beta^{-n} \int_{-\infty}^x e^{ru} dQ_n(u).$$

It is easy to verify that
$$G_n(x) = G^{n \star}(\sqrt{n}x), \hspace{0.05in} G(x) = \beta^{-1} \int_{-\infty}^x e^{su} dQ(u).$$

Assuming that random variable $Z$ obeys distribution function $G$, the following
relations are valid
\begin{equation} \label{eq:zrv}
 \Eb Z = \mu = \frac{m_1}{\beta}, \hspace{0.1in} \Eb (Z - \mu)^2 = \delta^2
 = \frac{m_2}{\beta} - \mu^2.
\end{equation}
Besides,
\begin{equation} \label{eq:ftail}
  1-Q_n(t) = \beta^n \int_t^{\infty} e^{-ru} dG_n(u), \hspace{0.1in}
  1 - \Phi(t) = \exp{\{\frac{r^2}{2}\}}\int_t^{\infty} e^{-ru} d \Phi(u-r).
\end{equation}

The following notations will be used below
\begin{equation} \label{eq:alfa}
   \alpha_k = t_0^{3-k} \tau^{-k} \exp{\{-0.5(1-\tau)^2t_0^2\}}, \hspace{0.1in} k=0..3;
\end{equation}
\begin{equation} \label{eq:alfa1}
  \log \Delta = \frac{c}{\tau^3 b^3}
  \exp{\{t_0^2 \left( \tau(1-\tau)+\frac{2(c-b)}{b^2} \right) \}};
\end{equation}
\begin{equation} \label{eq:alfa2}
   \mu=\frac{m_1}{\beta};
\end{equation}
\begin{equation} \label{eq:alfa3}
   \eta= |\mu| \sqrt{\overline{m}_2}(3\sqrt{\overline{m}_2} +|\mu| \sqrt{\beta}).
\end{equation}

\begin{lemma} \label{lm:one} The following bound is valid
\begin{equation} \label{eq:form5}
  \frac{\sqrt{n}}{\rho} | G_n(x\delta + \sqrt{n}\mu)-\Phi(x)| \leq 0.7655 q
\end{equation}
where
\begin{equation} \label{eq:form6}
q = \frac{1}{\beta\rho \delta^3} \int |y-\mu|^3 e^{sy} dQ(y) \leq
   \sqrt{\frac{\beta}{(\underline{m}_2-\beta \mu^2)^3}} \left(e^{sh}+\eta \right).
\end{equation}
\end{lemma}
\textit{Proof:} Inequality (\ref{eq:form5})
 follows from (\ref{eq:berry}), (\ref{eq:shig}) and (\ref{eq:zrv}).
Next, we consider upper bound (\ref{eq:form6})
$$q = \frac{1}{\beta\rho \delta^3} \int |y-\mu|^3 e^{sy} dQ(y) \leq
\frac{1}{\beta\rho \delta^3} \int (|y| + |\mu|) (y^2 - 2 \mu y + \mu^2) dQ(y) $$
$$\leq \frac{1}{\beta\rho \delta^3} \int
   \left( |y|^3+3|\mu|y^2+\mu^2|y| -|\mu|^3(2\beta-1)\right) e^{sy} dQ(y) $$
\begin{equation} \label{eq:form6a}
   \leq \frac{1}{\beta\rho \delta^3} \int
   \left( |y|^3+3|\mu|y^2+\mu^2|y| \right) e^{sy} dQ(y).
\end{equation}
The required bound may be deduced as a consequence of
(\ref{eq:zrv}) and (\ref{eq:form6a}) plus (\ref{eq:pow3}) and (\ref{eq:pow4}). $\blacksquare$

\begin{lemma} \label{lm:two} Suppose that
\begin{equation} \label{eq:term5}
   \beta \leq 1+\frac{\rho s^2h}{6}.
\end{equation}
Then,
\begin{equation} \label{eq:term2x}
 | \Phi(x\delta+\sqrt{n}\mu-r) - \Phi(x\delta) | \leq \frac{\rho}{h^2} e^{sh}.
\end{equation}
\end{lemma}
\textit{Proof:}
Based on the definition of normal distribution we have
\begin{equation} \label{eq:term2}
 | \Phi(x\delta+\sqrt{n}\mu-r) - \Phi(x\delta) | \leq
 \frac{|\sqrt{n}\mu - r|}{\sqrt{2\pi}} =\sqrt{\frac{n}{2\pi}} |\mu-s|.
\end{equation}

As far as $\beta \geq 1$, the following relation is valid according to (\ref{eq:pow1})
\begin{equation} \label{eq:term3}
  \mu-s \leq m_1 -s \leq \frac{\rho}{h^2} e^{sh}.
\end{equation}

The inequality $$m_1 \geq s- \frac{\rho}{h^2} \left( 1+sh + 0.5(sh)^2 \right)$$
follows from $x\exp{\{sx\}} \geq x(1+sx+0.5(sx)^2)$ $\forall x \in \mathbb{R}$,
and from the low bounds:
\begin{equation} \label{eq:lowbound}
    \Eb Y \geq -\frac{\rho}{h^2};
    \Eb Y^2 \geq 1-\frac{\rho}{h};
    \Eb Y^3 \geq -\rho.
\end{equation}
Therefore,
\begin{equation} \label{eq:term4}
   s-\mu \leq \frac{s-m_1}{\beta} +(\beta-1)s \leq \frac{\rho}{h^2}
   (1+sh+0.5(sh)^2 ) +(\beta-1)s \leq \frac{\rho}{h^2} e^{sh}
\end{equation}
subject to the condition (\ref{eq:term5}). $\blacksquare$

Next, we use the property $\rho \geq 1$, which follows from H\"older's inequality
applied to the Berry-Esseen condition $\Eb X^2 =1$.

In accordance with (\ref{eq:pow0}) and (\ref{eq:cond2})

$$\beta \leq 1+ \frac{t_0^2(1-\tau)^2}{2n} + \frac{\rho}{(\tau t_0 \sqrt{n})^3}
  \exp{\{ \tau (1-\tau) t_0^2 \}} $$
\begin{equation} \label{eq:term6}
   \leq 1+\frac{\gamma(t_0)}{2\sqrt{n}} \left(
   t_0^2(1-\tau)^2+\frac{2\gamma(t_0)\exp{\{\tau (1-\tau )t_0^2 \}}}
   {(\tau t_0)^3} \right).
\end{equation}

Combining (\ref{eq:term5}) and (\ref{eq:term6}) we obtain stronger condition
\begin{equation} \label{eq:term7}
   c \cdot \exp{\{\frac{2(c-b)t_0^2}{b^2} \}}
   \left( t_0^2+\frac{2c \cdot \exp{\{t_0^2(\tau (1-\tau) +\frac{2(c-b)}{b^2}) \}}}
   {b^3 \tau^3 (1-\tau)^2 } \right) \leq \frac{b^3 \tau}{3}.
\end{equation}

\begin{lemma} \label{lm:tri} The following upper bound is valid
\begin{equation} \label{eq:term9x}
  \underset{t \geq t_0}{\sup} | \Phi(t\delta)-\Phi(t) | \leq
  \frac{t_0 |\delta-1|
  \left[\exp{\{-\frac{t_0^2}{2}\}}+\exp{\{-\frac{(\delta t_0)^2}{2}\}}\right]}{2 \sqrt{2\pi}}.
\end{equation}
\end{lemma}
\textit{Proof:} The required inequality follows from convexity of the exponential function.

\begin{lemma} \label{th:second} (\texttt{Center approximation})
The following upper bound is valid
\begin{subequations}
\begin{align}
\label{eq:form7}
   |t|^3 H_n(t) \leq  B_C(t_0) = \tau^{-3} +\alpha_3 \Delta  \\
    \label{eq:form7xa}
   + \sqrt{\frac{2}{\pi}}
   \left[ \alpha_2 + 0.25 \cdot \alpha_1 t_0
   \left(\exp{\{-\frac{t_0^2}{2}\}}+ \exp{\{-\frac{(\delta t_0)^2}{2}\}} \right) \right]   \\
   \label{eq:form7x}
   + 1.531 \sqrt{\frac{\beta}{(\underline{m}_2-\beta\mu^2)^3}} \left( \alpha_0 +\eta t_0^3
   \exp{\{-\frac{t_0^2(1-\tau^2)}{2}\}} \right);
\end{align}
\end{subequations}

under conditions (\ref{eq:cond2}), (\ref{eq:term7}) and
\begin{subequations}
\begin{align}
  \label{eq:gcond}
  \frac{2(b-c)}{b^2} > \tau(1-\tau);   \\
  \label{eq:gcond6}
  t_0^2 \geq \left( \frac{2(b-c)}{b^2} -\tau(1-\tau)\right)^{-1};  \\
  \label{eq:gcond5}
  t_0^2 \geq  \frac{3b^2}{2(b-c)};             \\
  \label{eq:gcond1}
  t_0^2 \geq \frac{5}{2\tau(1-\tau)};         \\
  \label{eq:gcond7}
  t_0^2 \geq \frac{3}{(1-\tau)^2}
  \end{align}
\end{subequations}
where $0<\tau<1$ and $b>c \geq 1$.
\end{lemma}
\textit{Proof:} 
Using (\ref{eq:mark}) and (\ref{eq:cond2}) we obtain the upper bound
\begin{equation} \label{eq:bnd}
   | \beta -1 - \frac{r^2}{2n} | \leq \frac{\varepsilon}{n} \leq
   \frac{c}{n \tau^3 b^3} \exp{\{t_0^2 \left( \tau(1-\tau)+\frac{2(c-b)}{b^2} \right) \}}.
\end{equation}

Suppose that
\begin{subequations}
\begin{align}
 \label{eq:cond3}
    5 \varepsilon \leq r^2;   \\
 \label{eq:cond3a}
    r^2 \leq 2 \varepsilon \sqrt{n}.
\end{align}
\end{subequations}

Inequality (\ref{eq:cond3}) follows from stronger condition
\begin{equation} \label{eq:cond4}
    5 \cdot c \cdot \exp{\{t_0^2(\tau(1-\tau)+\frac{2(c-b)}{b^2})\}} \leq \tau^3(1-\tau)^2b^3t_0^2
\end{equation}
under (\ref{eq:gcond}): condition that the left
part of (\ref{eq:cond4}) is a non-increasing function of $t_0$
(means, the inequality will be valid $\forall t \geq t_0$).

Inequality (\ref{eq:cond3a}) follows from stronger condition
\begin{equation} \label{eq:cond6}
  (1-\tau)^2\tau^3t_0^5 \exp{\{-\tau(1-\tau)t_0^2\}} \leq 2
\end{equation}
under (\ref{eq:gcond1}): condition that the left
part of (\ref{eq:cond6}) is a non-increasing function of $t_0$.

According to (\ref{eq:pow0})  $\beta \leq \exp{\{\frac{s^2}{2} + \frac{\varepsilon}{n} \}}$.
Therefore,
\begin{equation} \label{eq:bnd1}
   \frac{\sqrt{n}}{\rho} e^{-rt} |\beta^n -e^{0.5r^2} | \leq
   \frac{\sqrt{n}}{\rho} \varepsilon \exp{\{\frac{r^2}{2}+\varepsilon -rt \}} \leq
   \alpha_3 \Delta t^{-3}
\end{equation}
where $\varepsilon \leq \log \Delta$ as an equivalent of (\ref{eq:cond2}),
and $\Delta$ is defined in (\ref{eq:alfa1}).

The inequality (\ref{eq:termvs}) is similar to (\ref{eq:markov}).
Then, we use representation (\ref{eq:ftail}) and bound (\ref{eq:bnd1})
\begin{equation} \label{eq:termvs}
  | H_n(t) | \leq \frac{\sqrt{n}}{\rho } |\Phi(t) - Q_n(t) | +\frac{n^{\frac{3}{2}}}{\rho }
  \mathbb{P}(|X|>h)
\end{equation}
$$\leq \frac{\sqrt{n}}{\rho} \left( |\beta^n - e^{0.5r^2} | \int_t^{\infty} e^{-ru} d G_n(u) +
  e^{0.5r^2} | \int_t^{\infty} e^{-ru} d \left[ G_n(u)-\Phi(u-r) \right] | \right) + (\tau t)^{-3}$$
$$\leq \frac{\sqrt{n}}{\rho} \left( | \beta^n - e^{0.5r^2} | e^{-rt} +
  2 \exp{\{\frac{r^2}{2}-rt\}} \sup_{x \in \mathbb{R}} |G_n(x)-\Phi(x-r)| \right) +
  (\tau t )^{-3}$$
\begin{equation} \label{eq:term}
   \leq \left( \alpha_3 \Delta +\tau^{-3} \right) t^{-3} + \frac{2\sqrt{n}}{\rho}
   \exp{\{\frac{r^2}{2}-rt\}} \sup_{x \in \mathbb{R}} |G_n(x) - \Phi(x-r) |.
\end{equation}

Consider the last term in (\ref{eq:term})
$$\sup_{x \in \mathbb{R}} | G_n(x) - \Phi(x-r) | \leq
 \sup_{x \in \mathbb{R}} \{ | \Phi(x\delta+\sqrt{n}\mu-r)-\Phi(x\delta)|
$$
\begin{equation} \label{eq:term1}
  + |\Phi(x\delta) - \Phi(x) | + |\Phi(x)-G_n(x\delta + \sqrt{n}\mu)| \}.
\end{equation}

Using (\ref{eq:term2x}) we deduce that
\begin{equation} \label{eq:term8}
   \frac{2\sqrt{n}}{\rho} \exp{\{ 0.5r^2-rt\}} |\Phi(x\delta+\sqrt{n}\mu-r)-
   \Phi(x\delta)| \leq \sqrt{\frac{2}{\pi}} \frac{\alpha_2}{t^3}
\end{equation}
under (\ref{eq:gcond7}): condition which ensure
that $\alpha_2$ is a non-increasing function of $t_0$.

Based on the general inequality
\[ |\delta -1| \leq 0.5 \max{\left(\delta^2-1, \frac{1-\delta^2}{\delta} \right)} \]
and result of the Lemma~\ref{lm:tri} we can conclude that
\begin{equation} \label{eq:term9}
   | \Phi(x\delta)-\Phi(x) | \leq \frac{ t_0 \zeta }{4\sqrt{2\pi}}
   \left[\exp{\{-\frac{t_0^2}{2}\}}+ \exp{\{-\frac{(\delta t_0)^2}{2}\}} \right]
\end{equation}
subject to the conditions
\begin{equation} \label{eq:cond7}
  \max{\{0.25, 1-0.5 \zeta \}} \leq \delta^2 \leq 1+\zeta,
  \zeta > 0.
\end{equation}
We have
\begin{equation} \label{eq:form}
  \underline{m}_2 =
  1-\frac{\rho}{h}(1+sh) \leq m_2 \leq 1+\frac{\rho}{h}e^{sh} =
  \overline{m}_2
\end{equation}
where left inequality is valid according to $e^{sx} \geq 1+sx$ $\forall x$,
and (\ref{eq:lowbound}).

The right inequality is valid according to (\ref{eq:pow2}).
It follows from (\ref{eq:form}) and
\[ \frac{1}{\beta} \geq 2-\beta, \beta \geq 1, \] that
\begin{equation} \label{eq:form1}
  \delta^2 \geq \frac{1}{\beta} -\frac{\rho}{h}(1+sh)-\mu^2 \geq
  2-\beta -\frac{\rho}{h}(1+sh)-\mu^2.
\end{equation}

Furthermore, by (\ref{eq:term3}) and (\ref{eq:term4}),
\begin{equation} \label{eq:form2}
  | \mu | \leq t_0 \cdot \gamma(t_0) \left( 1-\tau +
  \frac{c \cdot \exp{\{t_0^2(\tau(1-\tau) + \frac{2(c-b)}{b^2}) \}}}{\tau^2b^3}\right).
\end{equation}

By inserting \[ \zeta = \frac{\rho}{h}e^{sh} \] into (\ref{eq:term9})  we have
\begin{subequations}  \label{eq:form3}
\begin{align}
  \frac{2\sqrt{n}}{\rho} \exp{\{0.5r^2 -rt \}} \cdot
  |\Phi(x\delta)-\Phi(x)|  \\
    \label{eq:form3s}
\leq \frac{\alpha_1}{2 \sqrt{2 \pi} \cdot t^3} t_0
  \left[\exp{\{-\frac{t_0^2}{2}\}}+ \exp{\{-\frac{(\delta t_0)^2}{2}\}} \right]
\end{align}
\end{subequations}
under (\ref{eq:gcond7}): condition which ensure that that $\alpha_1$ is a
non-increasing function of $t_0$.

The following inequality was derived using (\ref{eq:cond2}) and (\ref{eq:form1})
and corresponds to the condition $\delta^2 \geq 0.25$ of (\ref{eq:cond7})
\begin{equation} \label{eq:form4}
  \beta-1+\mu^2+c \cdot t_0^2 b^{-3} \exp{\{2(c-b)\left(\frac{t_0}{b}\right)^2\}}
  \left[ \tau^{-1} +t_0^2(1-\tau)\right] \leq 0.75.
\end{equation}
under (\ref{eq:gcond1}): this condition (combined with condition (\ref{eq:gcond}))
will ensure that the left part of (\ref{eq:form4}) is a non-increasing.

Note that we can use weaker condition \[ t_0 \geq \frac{5b^2}{4(b-c)} \]
in order to ensure that the left part of (\ref{eq:form4}) is a non-increasing
as a function of $t_0$. But, we used already stronger
(according to (\ref{eq:gcond})) condition (\ref{eq:gcond1}),
which is an essential for (\ref{eq:cond6}). Respectively, we will leave in above
and further cases only condition, which is stronger.

Next inequality corresponds to $\delta^2 \geq 1-0.5 \cdot \zeta$ of (\ref{eq:cond7})
and was obtained using (\ref{eq:form1})
\begin{equation} \label{eq:form4bd}
   1+\frac{\rho}{h} \left( \frac{e^{sh}}{2} -1 -sh \right) \geq \beta+\mu^2.
\end{equation}
Then, we apply (\ref{eq:term6}) and (\ref{eq:form2}) to the right side of (\ref{eq:form4bd})
\begin{subequations}
\begin{align}
   \label{eq:form4ba}
   \tau \cdot t_0 \cdot \gamma(t_0) ( 0.5t_0^2(1-\tau)^2 + c  \cdot \frac{
  \exp{\{t_0^2(\tau(1-\tau)+\frac{2(c-b)}{b^2}) \}}}{b^3 \tau^3}   \\
   \label{eq:form4b}
   + t_0^2 \left( 1-\tau + c \cdot
   \frac{\exp{\{t_0^2( \tau(1-\tau)+\frac{2(c-b)}{b^2}) \}}}{b^3 \tau^2}
   \right)^2 ) \\
    \label{eq:form4c}
    \leq 0.5 \cdot \exp{\{t_0^2\tau(1-\tau)\}}-1-t_0^2\tau(1-\tau)
\end{align}
\end{subequations}
under condition (\ref{eq:gcond1}), which ensure that (\ref{eq:form4c}) is a non-decreasing;
plus (\ref{eq:gcond5}): condition that the left part represented by (\ref{eq:form4ba}) and
(\ref{eq:form4b}) is a non-increasing.
The condition $\delta^2 \leq 1+\zeta$ of (\ref{eq:term9}) is always valid according
to (\ref{eq:pow2}).

We can re-write estimators (\ref{eq:term6}) and (\ref{eq:form}) in a more detailed form
using ``center'' condition
(\ref{eq:cond2})
\begin{equation} \label{eq:form6b}
  \beta \leq 1+ \gamma^2(t_0) \left( \frac{(1-\tau)^2t_0^2}{2} +
  \frac{c \cdot
  \exp{\{t_0^2(\tau(1-\tau)+\frac{2(c-b)}{b^2}})\}}{\tau^3b^3}\right)
\end{equation}
under (\ref{eq:gcond}) and (\ref{eq:gcond1}): conditions, which ensure that the upper
estimator (\ref{eq:form6b}) is a non-increasing;
\begin{equation} \label{eq:form6c}
  \underline{m}_2 = 1-\frac{\gamma(t_0)}{t_0} (\frac{1}{\tau} +(1-\tau)t_0^2) \leq m_2 \leq
  1+ \frac{\gamma(t_0)\exp{\{t_0^2\tau(1-\tau)\}}}{t_0\tau} = \overline{m}_2
\end{equation}
under (\ref{eq:gcond6}): condition that the upper bound (\ref{eq:form6c}) is a non-increasing,
and under (\ref{eq:gcond1}): condition that the low bound (\ref{eq:form6c}) is a non-decreasing.

By (\ref{eq:form5}) and (\ref{eq:form6})
$$\frac{2\sqrt{n}}{\rho}\exp{\{0.5r^2-rt\}}|G_n(x\delta+\sqrt{n}\mu)-\Phi(x)| $$
\begin{equation} \label{eq:form8}
   \leq \frac{1.531}{t^3} \sqrt{\frac{\beta}{(\underline{m}_2-\beta\mu^2)^3}}
   \left(\alpha_0+\eta t_0^3\exp{\{-0.5t_0^2(1-\tau^2)\}}\right)
\end{equation}
under (\ref{eq:gcond7}): condition that $\alpha_0$ is a non-increasing.

Assuming that $c=1$, and combining (\ref{eq:term}), (\ref{eq:term1}), (\ref{eq:term8}),
(\ref{eq:form3}) and (\ref{eq:form8}) under conditions
(\ref{eq:gcond}) - (\ref{eq:gcond7}),
(\ref{eq:cond4}), (\ref{eq:cond6}),
(\ref{eq:term7}), (\ref{eq:form4}) and (\ref{eq:form4b}) we obtain required result. $\blacksquare$

\begin{proposition} \label{th:prop}
Suppose that the sample size $n$ is fixed. Then,
\begin{equation}
   \lim_{t \rightarrow \infty}  C(t) \leq 1
\end{equation}
under conditions of the Berry-Esseen Theorem (\ref{eq:beress}).
\end{proposition}

\textit{Proof} We need to consider (\ref{eq:fappr}) only if
sample size $n$ is fixed and $t$ is large enough.

Clearly, we can construct the functions $c(t)$ and $b(t)$: $1 < c(t) < b(t) \longrightarrow 1
 \hspace{0.1in} \text{if} \hspace{0.1in} t \rightarrow \infty,$
under ``tail'' condition (\ref{eq:frst})
$$b(t) \left( 1-\frac{b(t)}{2t^2}
  \log{ \frac{\sqrt{n}t^3}{\rho a} } \right) - c(t) \geq 0,$$
such that the upper bound for $t^3H_n(t)$ obtained from (\ref{eq:happr})
will tend to $1$ if $t \rightarrow \infty$.

\section{Concluding Remarks}  \label{sec:conclud}
According to \cite{Hal83}, the probability literature contains a large body of very elegant
mathematical theory which describes the rate of convergence in the central limit theorem.
These results often involve a uniform measure of the rate of convergence. Statisticians are
sometimes rather skeptical of such theory, pointing out that it is disjoint from the more
practical problems, which they encounter. Frequently, they are only interested in the rate
of convergence in isolated points.

For example, using the nonuniform bound (\ref{eq:nikul}) we can construct
in analytical form the upper bound for the
confidence interval based on the sample mean as an estimator of the location parameter
(\cite{Ben62} and \cite{HalJin95}):
$$\mathbb{P}(|\frac{1}{n} \sum_{i=1}^n X_i - \theta| > \varepsilon)
  \leq 2\left( 1 - \Phi(\sqrt{n}\varepsilon) +
  \frac{\rho C(\sqrt{n}\varepsilon)}{n^2\varepsilon^3} \right)$$
where $\sqrt{n}\varepsilon \geq 1$. The Table~\ref{tb:table1} demonstrates
advantage of the bound (\ref{eq:nikul}) if $\sqrt{n}\varepsilon \geq 3.3.$


\begin{thebibliography}{13}

\bibitem{Ben62}
G. Bennett. 
\lq\lq Probability Inequalities for the Sum of Independent Random Variables."
{\em Journal of the American Statistical Association}, vol. 57, pp. 33-45, 1962.

\bibitem{Ber41}
A. Berry.
\lq\lq The accuracy of the Gaussian approximation to the sum of independent variables."
{\em Thans. Amer. Math. Soc.}, vol. 49, pp. 122-136, 1941.

\bibitem{Ess45}
G. Esseen.
\lq\lq Fourier analysis of distribution functions. 
A mathematical study of the Laplace-Gaussian law."
{\em Acta Math.}, vol. 77, pp. 1-125, 1945.

\bibitem{Hal82}
P. Hall.
\lq\lq Rates of Convergence in the Central Limit Theorem."
{\em Pitman}, 1982.

\bibitem{Hal83}
P. Hall.
\lq\lq Sets which determine the rate of convergence in the central limit theorem." 
{\em The Annals of Probability}, vol. 11(2), pp. 355-361, 1983.

\bibitem{HalJin95}
P. Hall and B. Jing.
\lq\lq Uniform coverage bounds for confidence intervals and berry-esseen theorems
for edgeworth expansion."
{\em The Annals of Statistics}, vol. 23(2), pp. 363-375, 1995.

\bibitem{Mic76}
R. Michel.
\lq\lq Nonuniform Central Limit Bounds with Applications to Probabilities of Deviations."
{\em The Annals of Probability}, vol. 4(1), pp. 102-106, 1976.

\bibitem{Mic81}
R. Michel. 
\lq\lq On the Constant in the Nonuniform Version of the Berry-Essen Theorem."
{\em Z. Wahrscheinlichkcitstheorie verw. Geb.}, vol. 55, pp. 109-117, 1981.

\bibitem{Nag65}
S. Nagaev.
\lq\lq Some limit theorems for large deviations." 
{\em Theory of Probability and its Applications}, vol. 10, pp. 214-235, 1965.

\bibitem{Nag79}
S. Nagaev.
\lq\lq Large deviations of sums of independent random variables."
{\em The Annals of Probability}, vol. 7, pp. 745-789, 1979.

\bibitem{Nik92}
V. Nikulin.
\lq\lq Nonuniform bounds for the remainder term in the central limit theorem." 
{\em Theory of Probability and its Applications}, vol. 34(4), pp. 831-832, 1992.

\bibitem{NikPad98}
V. Nikulin and L. Paditz.
\lq\lq A note on nonuniform CLT-bound."
{\em 7th Vilnius Conference on Probability Theory}, pp. 358-359, 1998.

\bibitem{Shi82}
I. Shiganov.
\lq\lq Refinement of the upper bound of a constant in the remainder term of the
 central limit theorem."
{\em Stability Problems of Stochastic Models, Moscow},
Nauchno-Issled. Inst. Sistem. Issled., pp. 109-115, 1982.

\end{thebibliography}
\end{document}